\documentclass[12pt]{article}
\usepackage{amssymb}
\usepackage{amsmath,amsthm}
\usepackage[latin1]{inputenc}
\usepackage{epsfig}
\usepackage{psfig}
\usepackage{hyperref}

\usepackage{graphicx}

 \newtheorem{remark}{Remark}

 \newtheorem{lemma}[remark]{Lemma}
 \newtheorem{theorem}[remark]{Theorem}
 \newtheorem{proposition}[remark]{Proposition}
 \newtheorem{corollary}[remark]{Corollary}
 \newtheorem{example}[remark]{Example}

\title{The $(\alpha,\beta,s,t)$-diameter of  graphs:\\A particular case of conditional diameter}

\author{J. A. Rodr\'{\i}guez\footnote{e-mail:\mbox{\tt
    juanalberto.rodriguez\@@uc3m.es}}\\
{\em Departamento de Matem\'{a}ticas}\\
Universidad Carlos III de Madrid\\ Avda. de la Universidad 30,
28911 Leganés (Madrid),  Spain }

\date{}

\begin{document}

\maketitle

\begin{abstract}
The conditional diameter of a connected graph $\Gamma=(V,E)$ is
defined as follows: given a property ${\cal P}$ of a pair
$(\Gamma_1, \Gamma_2)$ of subgraphs of $\Gamma$, the so-called
\emph{conditional diameter}  or ${\cal P}$-{\em diameter} measures
the maximum distance among subgraphs satisfying ${\cal P}$. That
is,
\[
D_{{\cal P}}(\Gamma):=\max_{\Gamma_1, \Gamma_2\subset \Gamma} \{
\partial(\Gamma_1, \Gamma_2): \Gamma_1, \Gamma_2 \quad {\rm
satisfy }\quad {\cal P}\}.
\]
In this paper we consider the conditional diameter in which ${\cal
P}$ requires that $\delta(u)\ge \alpha$ for all $ u\in V(\Gamma_1)$,
$\delta(v)\ge \beta$ for all $v\in V(\Gamma_2)$, $| V(\Gamma_1)| \ge
s$ and $| V(\Gamma_2)| \ge  t$ for some integers $1\le s,t\le |V|$
and $\delta \le \alpha , \beta \le \Delta$, where $\delta(x)$
denotes the degree of a vertex $x$ of $\Gamma$, $\delta$ denotes the
minimum degree and $\Delta$  the maximum degree of $\Gamma$. The
conditional diameter obtained is called $(\alpha ,\beta
,s,t)$-\emph{diameter}. We obtain upper bounds on the $(\alpha
,\beta ,s,t)$-diameter by using the $k$-alternating polynomials  on
the mesh of eigenvalues of an associated weighted graph. The method
provides also bounds for other parameters such as vertex separators.
\end{abstract}



\section{Introduction}

In this paper all graphs $\Gamma=(V,E)$ will be finite, undirected,
simple and connected. The order of $\Gamma$, $|V(\Gamma )|$, will be
denoted by $n$, and the size, $|E(\Gamma )|$, will be denoted by
$m$. The degree of a vertex $v_i\in V(\Gamma)$ will be denoted by
$\delta(v_i)$ (or by $\delta_i$ for short), the minimum degree of
$\Gamma$ will be denoted by $\delta$ and the maximum by $\Delta$.
Moreover, the minimum degree of a vertex subset $U\subseteq
V(\Gamma)$ will be denoted by $\delta (U)$: $$ \delta(U):=\min_{u\in
U}\{\delta(u)\}.$$ We recall that the \emph{distance}
$\partial(u,v)$ between two vertices $u$ and $v$ is the mi\-nimum of
the lengths of paths between  $u$ and $v$ and the distance
$\partial(U, W)$ between two sets of vertices $U, W\subseteq
V(\Gamma)$ is defined as
$$\partial(U, W) := \min_{u\in U, v\in W} \lbrace \partial(u,v)\rbrace .$$

The conditional diameter of a graph was defined in \cite{BC} as
follows: given a property ${\cal P}$ of a pair $(\Gamma_1,
\Gamma_2)$ of subgraphs of $\Gamma$, the so-called
{\emph{conditional diameter}  or ${\cal P}$-{\em diameter}
measures the maximum distance among subgraphs satisfying ${\cal
P}$. That is,
\[
D_{{\cal P}}(\Gamma):=\max_{\Gamma_1, \Gamma_2\subset \Gamma} \{
\partial(V(\Gamma_1), V(\Gamma_2)): \Gamma_1, \Gamma_2 \quad {\rm
satisfy }\quad {\cal P}\}
\]
The study of conditional diameter is of interest, for instance, in
the design of interconnection networks when we need to minimize
the communication delays between the clusters represented by such
subgraphs. A direct application of conditional diameter to the
study of the superconnectivity of interconnection networks is
given in \cite{BC,degreDiamConn2} and \cite{degreDiamConn1} .

If ${\cal P}$ is the property of $\Gamma_i$, $i=1,2,$ being trivial
(that is, isolated vertices) the conditional diameter $D_{\cal
P}(\Gamma)$ coincides with the standard diameter $D(\Gamma)$.
Moreover, if ${\cal P}$ requires that $| V(\Gamma_1)| =s$ and $|
V(\Gamma_2)| = t$ for some integers $1\le s,t\le | V(\Gamma)|$, the
conditional diameter obtained is called $(s,t)$\emph{-diameter} and
denoted by $D_{(s,t)}(\Gamma)$. This conditional diameter  was
bounded by Fiol, Garriga and Yebra in \cite{dc} by using the
eigenvalues of the standard adjacency matrix.

In this paper we consider the case in which ${\cal P}$ requires that
$\delta (V(\Gamma_1))\ge \alpha$, $\delta(V(\Gamma_2))\ge\beta$, $|
V(\Gamma_1)| \ge s$ and $| V(\Gamma_2)| \ge t$ for some integers
$1\le s,t\le | V(\Gamma)|$ and $\delta \le \alpha , \beta \le
\Delta$. The conditional diameter obtained is called
$(\alpha,\beta,s,t)$-\emph{diameter} and will be denoted by
$D_{(s,t)}^{(\alpha,\beta)}(\Gamma)$. In particular, the
$(\alpha,\beta)$-\emph{degree} \emph{diameter} is defined by
\[
D^{(\alpha,\beta)}(\Gamma):=D_{(1,1)}^{(\alpha,\beta)}(\Gamma)=\max_{u,v\in
V}\{\partial(u,v):\delta(u)\ge \alpha, \delta(v)\ge \beta \}.
\]

In this paper we  obtain tight bounds on the
$(\alpha,\beta,s,t)$-diameter by using the $k$-alternating
polynomials  on the mesh of eigenvalues of a sui\-table adjacency
matrix that we call \emph{degree-adjacency matrix}. The method
provides also bounds for other parameters such as vertex separators.

\section{Degree-adjacency matrix}

 We define the $\emph{degree-adjacency matrix}$ of a
graph $\Gamma$ of order $n$ as the $n\times n$ matrix ${\cal A}$
whose ($i,j$)-entry is
\[ a_{ij}= \left\lbrace \begin{array}{ll} \frac{1}{\sqrt{\delta_i\delta_j}}  & {\rm if }\quad  v_i\sim v_j ;
                                                    \\
                                                0 &  {\rm otherwise} \end{array}
                                                \right. \]

The matrix ${\cal A}$ can be regarded as the adjacency matrix of a
weighted graph in which the edge-weight of the edge $v_iv_j$ is
equal to $\frac{1}{\sqrt{\delta_i\delta_j}}$, thus justifying the
terminology used. The degree-adjacency matrix is the adjacency
matrix derived from the Laplacian matrix used systematically by Fan
R. K. Chung \cite{Chung}.

In the case of $\nu=(\sqrt{\delta_1}, \sqrt{\delta_2}, ...,
\sqrt{\delta_n})$  we have ${\cal A} \nu =\nu$. Thus, $\lambda=1$
is an eigenvalue of ${\cal A}$ and $\nu$ is an eigenvector
associated to $\lambda$.  Hence, as ${\cal A}$ is non-negative and
irreducible in the case of connected graphs, by the
Perron-Frobenius theorem,  $\lambda=1$ is a simple eigenvalue and
$\lambda=1\ge | \lambda_j |$ for every eigenvalue $\lambda_j$ of
${\cal A}$. Hereafter the eigenvalues of ${\cal A}$ will be called
{\em degree-adjacency eigenvalues of} $\Gamma$.

It is well-known that there are non-isomorphic graphs that have the
same standard adjacency eigenvalues with the same multiplicities
(the so called cospectral graphs). For instance, two connected
graphs, both having the characteristic polynomial
$P(x)=x^6-7x^4-4x^3+7x^2+4x-1$, are shown in Figure \ref{ej1}. So,
in such cases, the spectral study doesn't allow to obtain structural
properties that differentiate both graphs.  Therefore, we can try to
study cospectral graphs by using an alternative matrix, for
instance, the degree-adjacency matrix ${\cal A}$. If we consider the
matrix ${\cal A}$, as one might expect, the eigenvalues of both
graphs are different: the left hand side graph has degree-adjacency
eigenvalues 1, $\pm\frac{1}{2}$ and $-\frac{1}{4}\left(1\pm
\sqrt{2.6} \right)$ (where the eigenvalue $-\frac{1}{2}$ has
multiplicity 2), on the other hand, the right hand side graph has
degree-adjacency eigenvalues 1, $\frac{-1\pm \sqrt{2}}{3}$, $\pm
\frac{\sqrt{3}}{3}$ and $- \frac{1}{3}$. Even so, the
degree-adjacency eigenvalues do not determine the graph. That is,
there are non-isomorphic graphs (and non-cospectral) that are
cospectral with regard to the degree-adjacency matrix. For instance,
the degree-adjacency eigenvalues of the cycle graph $C_4$ and the
semi-regular bipartite graph $K_{1,3}$ are the same: $1,0,0,-1$.
However, the standard eigenvalues  are $2,0,0,-2$, in the case of
$C_4$, and $\sqrt{3},0,0,-\sqrt{3}$ in the case of $K_{1,3}$.

\begin{figure}[h]

\begin{center}
\caption{Two cospectral graphs but not cospectral with regard to
${\cal A}$} \label{ej1}
\vspace{0,5cm}
\includegraphics[angle=90, width=7.2cm]{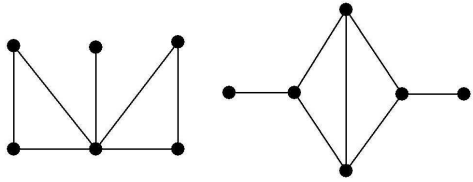}
\end{center}
\vspace*{-0.8cm}
\end{figure}

We identify  the degree-adjacency matrix ${\cal A}$ with an
endomorphism of the ``vertex-space" of $\Gamma$, $l^2(V(\Gamma))$
which, for any given indexing of the vertices, is isomorphic to
$\mathbb{R}^n$. Thus, for any vertex $v_i \in V(\Gamma)$, $e_i$
will denote the corresponding  unit vector of the canonical base
of $\mathbb{R}^n$.

If for two vertices $v_i,v_j\in V(\Gamma)$ we have $\partial
(v_i,v_j)>k$ then $({\cal A}^k(\Gamma))_{ij}=0$. Thus, for a real
polynomial $P$ of degree $k$, we have
\begin{equation}
\partial (v_i,v_j)>k \Rightarrow  P({\cal A}(\Gamma))_{ij}=0 .\label{just}
\end{equation}
Through this fact  we will study the $(\alpha,\beta,s,t)$-diameter
of $\Gamma$ by using the degree-adjacency matrix (or its
eigenvalues) and the $k$-alternating polynomials.

Another application of the degree-adjacency matrix can be found in
\cite{Randic} where spectral-like bounds on the higher Randi\'{c}
index $R_2(\Gamma)$ are given.

\section{Alternating polynomials}

The k-alternating polynomials, defined and studied in
\cite{alternante} by Fiol, Garriga and Yebra, can be defined as
follows:
 let ${\cal M}= \{\mu_1> \cdots >\mu_b\}$ be a mesh of real numbers. For any
 $k=0,1,...,b-1$ let $P_k$ denote the {\em k-alternating polynomial} associated to  ${\cal M}$.
 That is, the polynomial of  $\mathbb{R}_k[x]$
such that
 $$ P_k(\mu)= \sup_{P\in \mathbb{R}_k[x]} \left\{ P(\mu):  \| P \|_{\infty} \le 1 \right\} $$
 where  $\mu$ is any real number greater than $\mu_1$ and  $\| P \|_{\infty} =
 {\displaystyle\max_{1 \le i \le b}}\{| P(\mu_i)| \}$. We collect here some of its main properties, referring
 the reader to \cite{alternante} for a more detailed study.
  \begin{itemize}
  \item{  For any $k=0,1,...,b-1$ there is a unique $P_k$ which, moreover, is independent of the value of
  $\mu (> \mu_1)$;}
 \item{ $P_k$ has degree k;}
 \item{ $P_0(\mu)=1< P_1(\mu)< \cdots < P_{b-1}(\mu)$;}
 \item{ $P_k$ takes $k+1$  alternating values $\pm 1$ at the mesh
 points;}

\item{ There are explicit formulae for $P_0(=1),$ $P_1,$ $P_2,$ and $P_{b-1},$ while the other
polynomials can be computed by solving a linear programming
problem (for instance by the simplex method).}
 \end{itemize}

Hereafter the different eigenvalues of ${\cal A}$ will be denoted by
$\lambda_0=1,\lambda_1,\cdots,\lambda_b$ with
$\lambda_0=1>\lambda_1>\cdots>\lambda_b$.

\begin{proposition} \label{prop}
 Let $\Gamma$ be a simple and connected graph. Let  $P_k$ be the $k$-alternating polynomial associated
 to the mesh ${\cal M}= \{\lambda_1> \cdots >\lambda_b\}$ of degree-adjacency eigenvalues of $\Gamma$.
 Let $\nu$ be an eigenvector belonging to the eigenvalue $\lambda_0=1$.
 If $z\in \nu^\perp$ then $  \| P_k({\cal A}(\Gamma))z\|
\le \| P_k\|_{\infty} \| z
 \|$.
\end{proposition}

\begin{proof}
  Using the following decomposition of
 the vector $z$
$$z =\sum_{l=1}^b z_{l}, \quad \mbox{were} \quad  z_{l}\in \mbox{\rm Ker}({\cal A} - \lambda_l{\bf I}). $$
 we obtain
\begin{align*}
 \left\| P_k({\cal A})z\right\|^2 &= \left\| P_k({\cal A})\sum_{l=1}^b
 z_{l}\right\|^2
=\left\| \sum_{l=1}^b P_k(\lambda_l)z_{l}\right\|^2\\
&=  \sum_{l=1}^b (P_k(\lambda_l))^2 \| z_{l}\|^2 \le \|
P_k\|_{\infty}^2\sum_{l=1}^b  \| z_{l}\|^2
 = \| P_k\|_{\infty}^2 \| z\|^2.
 \end{align*}
 Hence, the result follows.
\end{proof}

Recently, the $k$-alternating polynomials have been successfully
applied to the study of several parameter related to the concept
of distance in graphs and hypergraphs. For instance, we cite
\cite{alternante,dc,dydm,cut,exceso,hyp}. We emphasize the
following result on the $(s,t)$-diameter \cite{dc}
\begin{equation} \label{cotaCD}
P_k(\lambda)>\sqrt{\left(\frac{\| v\|^2}{s}-1\right)\left(\frac{\|{
v}\|^2}{t}-1\right)}
  \Rightarrow D_{(s,t)}(\Gamma) \le k,
  \end{equation}
  where $P_k$ denotes the $k$-alternating polynomials on the mesh
  of eigenvalues of the standard adjacency matrix of
  $\Gamma$, $\lambda$ denotes the largest eigenvalue of $\Gamma$, and ${ v}$
  the eigenvector associated to $\lambda$ with minimum component
  1. In the case of regular graphs, as { v}={ j}, the all-1
  vector, the result (\ref{cotaCD}) simplifies to
\begin{equation} \label{cotaCDreg}
P_k(\lambda)>\sqrt{\left(\frac{n}{s}-1\right)\left(\frac{n}{t}-1\right)}
  \Rightarrow D_{(s,t)}(\Gamma) \le k,
  \end{equation}
where $n=| V(\Gamma)|.$

\section{Bounding the $(\alpha,\beta,s,t)$-diameter}
\begin{lemma} \label{lema}
Let $\Gamma$ be a simple and connected graph of size $m$. Let $P_k$
be the $k$-alternating polynomial associated to the mesh ${\cal M}=
\{\lambda_1> \cdots >\lambda_b\}$ of  degree-adjacency eigenvalues
of $\Gamma$.  Let $S=\{v_{i_1},v_{i_2},...,v_{i_s}\}$ and
$T=\{v_{j_1},v_{j_2},...,v_{j_t}\}$ be two sets of vertices of
$\Gamma$, and let
$\rho_s=\left(\sum_{l=1}^{s}\sqrt{\delta(v_{i_l})}\right)^2$,
$\rho_t=\left(\sum_{r=1}^{t}\sqrt{\delta(v_{j_r})}\right)^2$. Then,
$$P_k(1) >\sqrt{\left(\frac{2ms}{\rho_s}-1\right)\left(\frac{2mt}{\rho_t}-1\right)}
  \Rightarrow \partial(S,T) \le k.$$
\end{lemma}

\begin{proof}
 Let $\sigma=\sum_{l=1}^s e_{i_l}$ and $\tau=\sum_{r=1}^t e_{j_r}$ be
the vectors of $\mathbb{R}^n$ associated to the sets $S$ and $T$.
Using the following decomposition
\begin{equation}\label{desc}
\sigma=\frac{\langle \sigma,\nu\rangle}{\| \nu\|^2}\nu+
u=\frac{\sqrt{\rho_s}}{2m}\nu+ u, \quad \tau=\frac{\langle
\tau,\nu\rangle}{\| \nu\|^2}\nu+ w=\frac{\sqrt{\rho_t}}{2m}\nu+ w,
\end{equation}
where $\nu=(\sqrt{\delta_1}, \sqrt{\delta_2}, ...,
\sqrt{\delta_n})$  and $u,w\in\nu^{\perp},$ we obtain
\begin{align*}
\partial(S,T) >k
  &\Rightarrow \langle P_k({\cal A}) \sigma, \tau \rangle =0 \\
  &\Rightarrow P_k(1)\frac{\sqrt{\rho_s \rho_t}}{2m}=-\langle P_k ({\cal A})u, w \rangle.
\end{align*}
Thus, by the Cauchy-Schwarz inequality we have
$$\partial(S,T) >k \Rightarrow P_k(1)\frac{\sqrt{\rho_s\rho_t}}{2m} \le \| P_k({\cal A}) u\| \| w
\|,$$ and by Proposition \ref{prop} we obtain
\begin{equation}
\partial(S,T)
>k\Rightarrow P_k(1)\frac{\sqrt{\rho_s\rho_t}}{2m} \le
\|P_k\|_{\infty} \|u\| \|w\|.\label{ant}
\end{equation}
 Moreover, the
decomposition (\ref{desc}) leads to
$$  s=\| \sigma\|^2 = \frac{\rho_s}{2m}+ \| u\|^2\Rightarrow \| u\|=\sqrt{s-\frac{\rho_s}{2m}}$$
and
$$  t=\| \tau\|^2 = \frac{\rho_t}{2m}+ \| w\|^2\Rightarrow \| w\|=\sqrt{t-\frac{\rho_t}{2m}}$$
So, by (\ref{ant}), we obtain
\begin{equation}\label{final}
\partial(S,T) >k \Rightarrow  P_k(1)\sqrt{\rho_s\rho_t}\le
\sqrt{(2ms-\rho_s)(2mt-\rho_t)}.
\end{equation}
The converse of (\ref{final}) leads to the result.
\end{proof}

\begin{theorem} \label{teorema}
Let $\Gamma$ be a simple and connected graph of size $m$. Let $P_k$
be the $k$-alternating polynomial associated to the mesh ${\cal M}=
\{\lambda_1> \cdots >\lambda_b\}$ of degree-adjacency eigenvalues of
$\Gamma$. Then,
\begin{equation} \label{cotaFund}
P_k(1)>\sqrt{\left(\frac{2m}{s\alpha}-1\right)\left(\frac{2m}{t\beta}-1\right)}
  \Rightarrow D^{(\alpha,\beta)}_{(s,t)} \le k.
  \end{equation}
\end{theorem}

\begin{proof}
Consider $S,T\subset V(\Gamma)$ such that $|S|\ge s$,
$\delta(S)\ge \alpha$, $|T|\ge t$ and $\delta(T) \ge \beta$. Then
we have
$$\rho_s=\left(\sum_{u\in S}\sqrt{\delta(u)}\right)^2\ge s^2\alpha
\quad{\rm and} \quad \rho_t=\left(\sum_{v\in
T}\sqrt{\delta(v)}\right)^2\ge t^2\beta.$$ Hence,
$$\left(\frac{2m}{s\alpha}-1\right)\left(\frac{2m}{t\beta}-1\right)\ge
\left(\frac{2ms}{\rho_s}-1\right)\left(\frac{2mt}{\rho_t}-1\right).$$
Therefore, by Lemma \ref{lema}  the result follows.
\end{proof}

As we can see in the following examples, the above bound is attained
for several values of the related parameters.

\begin{figure}[h]
\begin{center}
\caption{ } \label{ej} \vspace{-1cm}
\includegraphics[angle=90, width=8cm]{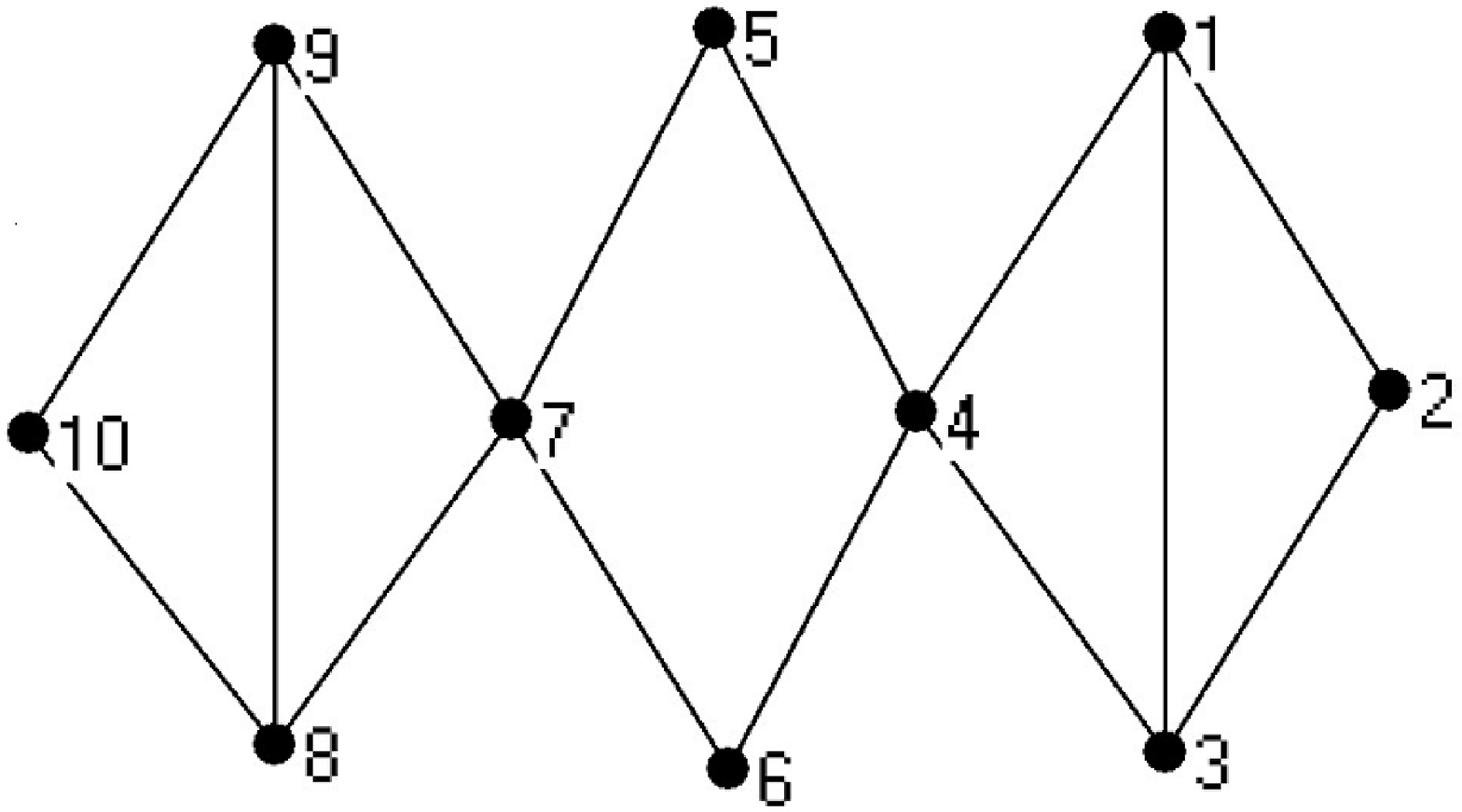}
\end{center}
\vspace*{-2cm}
\end{figure}

\begin{example}{\rm
The graph of Figure \ref{ej} has degree-adjacency eigenvalues:
\begin{center}
$\left\{1,
\frac{1+\sqrt{19}}{6},0.5358...,0,0,-\frac{1}{3},-\frac{1}{3},-0.3765...,\frac{1-\sqrt{19}}{6},-0.8259...\right\}$
\end{center}
from which we obtain $P_4(1)=3,89$, $P_5(1)=12,2$ and
$P_6(1)=266,5$. Thus, the following bounds are attained:
$D(\Gamma)=D_{(1,1)}^{(2,2)}(\Gamma)=D^{(2,2)}(\Gamma)\le 6$,
$D_{(1,1)}^{(2,3)}(\Gamma)=D^{(2,3)}(\Gamma)\le 5$,
$D_{(1,2)}^{(2,3)}(\Gamma)\le 5$, $D_{(2,2)}^{(3,3)}(\Gamma)\le 4$
and $D_{(3,3)}^{(2,2)}(\Gamma)=D_{(3,3)}(\Gamma)\le 4$. }
\end{example}

\begin{figure}[h]
\begin{center}
\caption{ } \label{ejemplo} \vspace{0,1cm}
\includegraphics[angle=90, width=7cm]{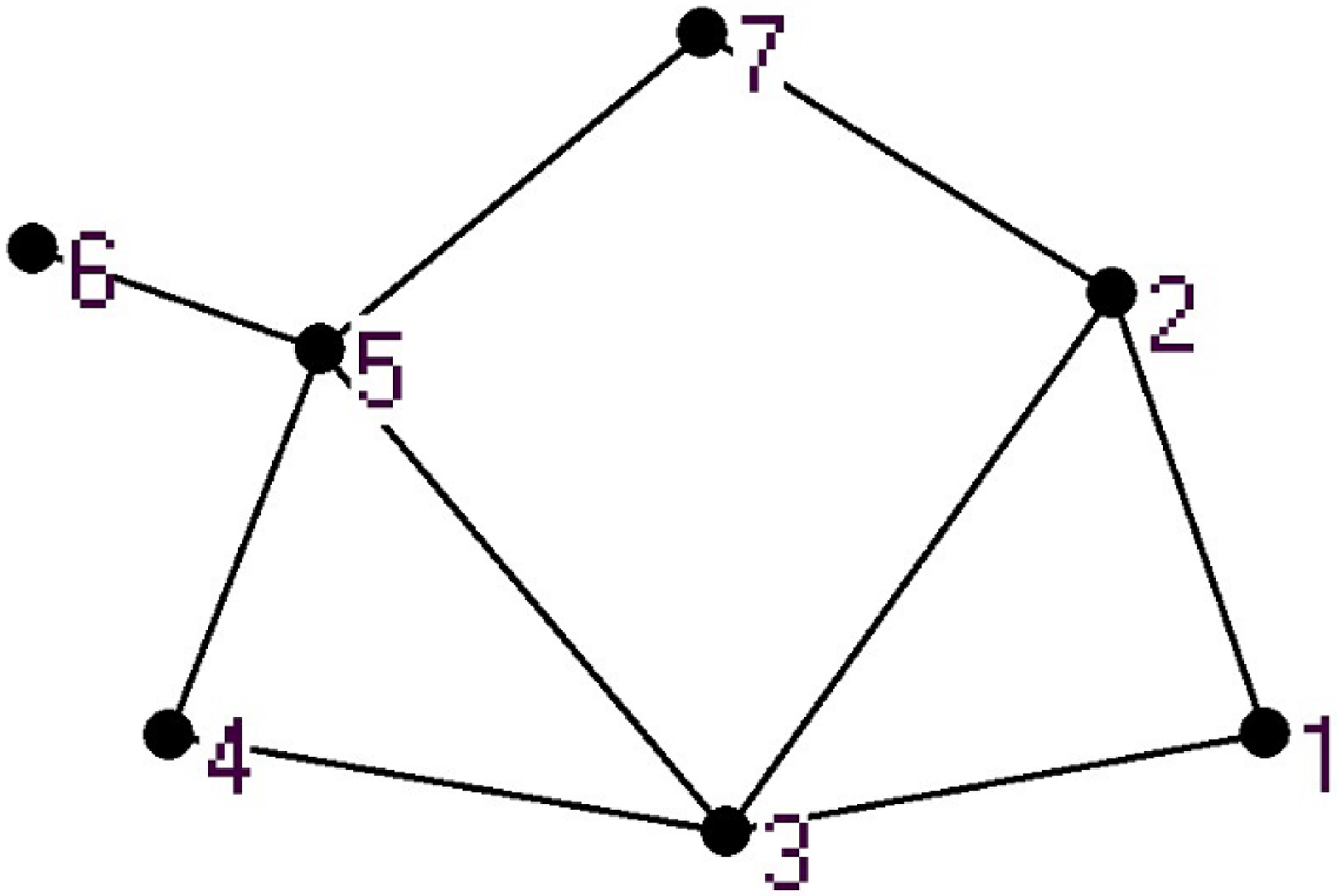}
\end{center}
\end{figure}

\begin{example}
{\rm
 The graph of Figure \ref{ejemplo} has degree-adjacency eigenvalues:

\begin{center}
$\left\{1,
\frac{-3+\sqrt{249}}{24},\frac{1}{4},0,-\frac{1}{2},-\frac{1}{2},\frac{-3-\sqrt{249}}{24}\right\}$
\end{center}
from which we obtain $P_1(1)=1.7$, $P_2(1)=5$, $P_3(1)=15.2$ and
$P_4(1)=58$. Thus, the following bounds are attained:
\begin{center}
 $D^{(1,2)}(\Gamma)\le 3$, $D^{(3,4)}(\Gamma) \le 2$ and $D^{(4,4)}(\Gamma) \le 1$.
\end{center}
}
\end{example}

 As particular cases of above theorem we derive the
following results in which the expression (\ref{cotaFund}) is
simplified.

\begin{corollary}
Let $\Gamma=(V,E)$ be a simple and connected graph of order $n$ and
size $m$. Let $P_k$ be the $k$-alternating polynomial associated to
the mesh ${\cal M}= \{\lambda_1> \cdots >\lambda_b\}$ of
degree-adjacency eigenvalues of $\Gamma$. Then,

\begin{itemize}
\item[{\rm (a)}]{
$P_k(1)>\frac{2m}{\alpha}-1 \Rightarrow
D^{(\alpha,\alpha)}(\Gamma) \le k. $}

\item[{\rm (b)}]{ The standard diameter is bounded by \\
$P_k(1)>\frac{2m}{\delta}-1 \Rightarrow D(\Gamma) \le k.$ }
\item[{\rm (c)}]{If $\Gamma$ is regular, the standard diameter is bounded
by \\ $P_k(1) >n-1 \Rightarrow D(\Gamma) \le k. $}
\item[{\rm (d)}]{If $\Gamma$ is an unicyclic graph, i.e., a connected graph containing exactly one cycle,
the standard diameter is bounded by \\ $P_k(1) >2n-1 \Rightarrow
D(\Gamma) \le k. $}
\item[{\rm (e)}]{If $\Gamma$ is regular, the
$(s,t)$-diameter is bounded by
$$P_k(1)>\sqrt{\left(\frac{n}{s}-1\right)\left(\frac{n}{t}-1\right)}
  \Rightarrow D_{(s,t)}(\Gamma) \le k.$$}
\end{itemize}
\end{corollary}

The bound (c) is an analogous result to the previous one given by
Fiol, Garriga and Yebra in \cite{alternante} by using the standard
adjacency matrix. Moreover, bound (e) is an analogous result to
(\ref{cotaCDreg}).

\subsection{Cutsets}

Now we are going to give some other consequences of above study
involving sets of vertices of equal cardinality and cut sets.

\begin{proposition} \label{PropIgualCardinal}
Let $\Gamma=(V,E)$ be a simple and connected graph of size $m$. Let
$P_k$ be the $k$-alternating polynomial associated to the mesh
${\cal M}= \{\lambda_1> \cdots >\lambda_b\}$ of degree-adjacency
eigenvalues of $\Gamma$. Let $S_1, S_2 \subset V(\Gamma)$ such that
$\vert S_1\vert = \vert S_2\vert=s$, $\partial (S_1, S_2)>k$ and
$\delta(v)\ge\alpha$ for all $v\in S_1 \cup S_2$. Then
\begin{equation} \label{igualcardinal}
   s \le \left\lfloor \frac{2m}{\alpha(P_k(1) +
1)}\right\rfloor \cdot
\end{equation}

\end{proposition}

\begin{proof}
Taking $\alpha=\beta$ and $s=t$, the converse of (\ref{cotaFund})
gives
$$\partial(S_1,S_2)> k\Rightarrow P_k(1)\le \frac{2m}{s\alpha}-1.$$
Solving for $s$, and considering that it is an integer, we obtain
the result.
\end{proof}

\begin{example}
{\rm  To show the tightness of above bound we consider again the
graph of  Figure \ref{ej}. For instance, taking $s=\vert S_1\vert
= \vert S_2\vert$, $\partial (S_1, S_2)>5$ and $\delta(v)\ge 2$
for all $v\in S_1 \cup S_2$, we obtain $s\le 1.$ Moreover, as for
this graph $P_3(1)=2,33$,  in the case of $\partial (S_1, S_2)>3$
and $\delta(v)\ge 3$ for all $v\in S_1 \cup S_2$, we obtain $s\le
2$. }
\end{example}

Note that, as in Proposition \ref{PropIgualCardinal}, if there are
two sets  $S_1, S_2 \subset V(\Gamma)$ such that $\vert S_1\vert =
\vert S_2\vert=s$, $\partial (S_1, S_2)>k$ and
$\delta(v)\ge\alpha$ for all $v\in S_1 \cup S_2$, then
\begin{equation} \label{gradoigualcardinal}
   \alpha \le \left\lfloor \frac{2m}{s(P_k(1) +
1)}\right\rfloor \cdot
\end{equation}

In the case of regular graphs, Proposition \ref{PropIgualCardinal}
 allows us to derive the following result.

 \begin{corollary} Let $\Gamma$ be
 a regular graph of order $n$.   Let $S_1, S_2
\subset V(\Gamma)$ such that $\vert S_1\vert = \vert S_2\vert=s$
and $\partial (S_1, S_2)>k$. Then
\begin{equation} \label{IgualCardRegular}
  s \le \left\lfloor \frac{n}{P_k(1) +
1}\right\rfloor \cdot
\end{equation}
\end{corollary}
The above result is analogous  to the previous one  given by Yebra
and the author in \cite{exceso}, for not necessarily regular
graphs, by using the standard Laplacian matrix. These result
becomes the main tool to the study of cut sets in \cite{cut}.

A {\em $k$-vertex separator} is a subset  $T_k \subset V(\Gamma)$
  whose deletion separates $V(\Gamma)$  into
 two sets of equal cardinality,  that are at distance greater than $k$.
We denote by $vs_k(\Gamma)$ the minimum cardinality among all
$k$-vertex  separators, that is,
 $$vs_k(\Gamma) = \min \{ \vert T_k\vert : T_k \mbox{\rm \mbox{ } is a $k$-vertex separator of } \Gamma \}.$$
 In \cite{cut} were obtained bounds on  $vs_k(\Gamma)$ by using  the $k$-alternating
 polynomials and the standard Laplacian spectrum. Proposition \ref{PropIgualCardinal} allows us to study
 a particular case of vertex separator: a $(\alpha,k)$-\emph{vertex separator}  is a vertex set
whose deletion separates $V(\Gamma)$  into two sets, $U$ and $W$,
of equal cardinality whose minimum vertex degree is $\alpha$, such
that $\partial(U,W)>k$. We denote by $vs_{(\alpha,k)}(\Gamma)$,
the minimum cardinality among all $(\alpha,k)$-vertex  separators.

\begin{corollary} \label{vertexSeparator}
Let $\Gamma=(V,E)$ be a simple and connected graph of order $n$ and
size $m$. Let $P_k$ be the $k$-alternating polynomial associated to
the mesh ${\cal M}= \{\lambda_1> \cdots >\lambda_b\}$ of
degree-adjacency eigenvalues of $\Gamma$. Then
\begin{equation} \label{cut}
vs_{(\alpha,k)}(\Gamma)\ge n-2\left\lfloor\frac{2m}{\alpha(P_k(1)
+ 1)}\right\rfloor \cdot
\end{equation}
\end{corollary}

\section{Laplacian matrix}

Now we consider the Laplacian matrix, $\emph{L}$, defined by Fan R.
K. Chung as $\emph{L}=I-{\cal A}$, where $I$ denotes the identity
matrix.

We denote by $\mu_0=0<\mu_1<\cdots<\mu_b$ the different eigenvalues
of \emph{L}. Thus, the eigenvalues of both matrices, $\emph{L}$ and
${\cal A}$, are related by
$$\mu_l=1-\lambda_l, \quad l=0,1,...,b.$$
 Notice also that the eigenvalue $\mu_0 = 0$ has  eigenvector
$\nu=(\sqrt{\delta_1},...,\sqrt{\delta_n})$ and multiplicity one
in the  case of connected graphs. Hence, both matrices, ${\cal A}$
and \emph{L}, lead to equivalent spectral-like results.
Particularly, the following theorem is the analogous of Theorem
\ref{cotaFund}. The proof is basically as before.

\begin{theorem}
Let $\Gamma=(V,E)$ be a simple and connected graph of size $m$. Let
$P_k$ be the $k$-alternating polynomial associated to the mesh
${\cal M}= \{\mu_1< \cdots <\mu_b\}$ of $\emph{L}=\emph{L}(\Gamma)$.
Then,
$$P_k(0) >\sqrt{\left(\frac{2m}{s\alpha}-1\right)\left(\frac{2m}{t\beta}-1\right)}
  \Rightarrow D^{(\alpha,\beta)}_{(s,t)}(\Gamma) \le k.$$
\end{theorem}

We recall that if we use the standard adjacency matrix and the
standard Laplacian matrix, the results are equivalent only in the
regular case. In this sense, a comparative study was done in
\cite{dydm}.

\end{document}